\documentclass[12pt, leqno]{article}
\usepackage{amsfonts}
\usepackage{amsmath}
\usepackage{amssymb}
\usepackage{amsthm}
\usepackage{amscd}
\usepackage{mathrsfs}
\usepackage{color}
\usepackage{latexsym}
\usepackage[dvipdfmx]{graphicx}
\usepackage{tikz}
\input cyracc.def
\newfont{\tencyr}{wncyr10}
\usepackage{pgf,pgfplots}

\usetikzlibrary{calc} 
\usetikzlibrary{intersections,calc,arrows.meta}
\usetikzlibrary{hobby}
\usetikzlibrary{decorations.markings}
\usetikzlibrary{knots}
\usetikzlibrary{3d}
\tikzset{otmm/.style={x={(-135:0.5cm)},y={(0:1cm)},z={(90:1cm)}}}
\usetikzlibrary{math}

\begin{document}

\begin{center}
{\Large \bf  Regularized determinant formulas for the zeta functions of 3-dimensional Riemannian foliated dynamical systems}
\end{center}

\vspace{.4cm}

\begin{center}
 Jes\'{u}s A. \'{A}lvarez L\'{o}pez, Junhyeong Kim, Masanori Morishita
\end{center}

\footnote[0]{2010 Mathematics Subject Classification: 37C25, 37C30, 53C12.\\
Key words: Riemannian foliated dynamical system, reduced leafwise cohomology,  dynamical zeta function, dynamical Lefschetz trace formula, regularized determinant}
\vspace{0.08cm}
\\
{\small
{\bf Abstract}: We prove a regularized determinant formula for the zeta functions of certain 3-dimensional Riemannian foliated dynamical systems, in terms of the infinitesimal operator induced by the flow acting on the reduced leafwise cohomologies. It is the formula conjectured by Deninger. The proof is based on relating dynamical  spectral $\xi$-function, analogues of the $\xi$-function in analytic number theory,  with the zeta function, by applying the distributional dynamical Lefschetz trace formula. 
}

\vspace{0.8cm}

\begin{center}
{\bf Introduction}  
\end{center}

In this paper, we prove a regularized determinant formula for the zeta functions of certain 3-dimensional Riemannian foliated dymanical systems. For this we employ the reduced leafwise cohomology and the infinitesimal operator induced by the flow. This type of formula was conjectured by Deninger ([D11], [D12]). It was motivated by the cohomological study of the zeta functions in number theory and  arithmetic geometry, which we recall in the following.

It is well known that \'{e}tale or crystalline cohomology is a powerful tool for the study of the zeta functions of schemes. For a basic example, let $X$ be a smooth proper curve $X$ over a finite field $\mathbb{F}_q$. The zeta function of $X$ is then defined by the Euler product
$$ \zeta_{X}(s) = \prod_{x \in |X|} (1 - {\rm N}(x)^{-s})^{-1}$$
over the set $|X|$ of closed points of $X$, where ${\rm N}(x) = \#\kappa(x)$, $\kappa(x)$ being the finite residue field of $x$, and the infinite product converges absolutely for ${\rm Re}(s) > 1$. Grothendieck associated to $X$ the \'{e}tale cohomology groups $H^i_{\footnotesize \mbox{\'{e}t}}(\overline{X}, \mathbb{Q}_{l})$ ($\overline{X} = X \times_{\mathbb{F}_q} \overline{\mathbb{F}}_q$, $\mathbb{Q}_l$ being the $l$-adic field with $l \nmid q$ and $0 \leq i \leq 2$), which are finite dimensional  $\mathbb{Q}_l$-vector spaces equipped with the action of the $q$-th power Frobenius automorphism ${\rm Fr}$, and he showed the following determinant formula 
$$ \zeta_X(s) = \prod_{n=0}^2 \det ({\rm id} -  {\rm Fr} \, q^{-s} \; | \; H^n_{\footnotesize \mbox{\'{e}t}}(\overline{X}, \mathbb{Q}_l))^{(-1)^{n+1}} \leqno{(0.1)}$$ 
in [G]. The key ingredient to prove (0.1) is Lefschetz trace formula: for any $k \in \mathbb{N}$, 
$$ \sum_{n=0}^2 (-1)^i {\rm tr}({\rm Fr}^k \; | \; H^n_{\footnotesize \mbox{\'{e}t}}(\overline{X}, \mathbb{Q}_l)) = \# X(\mathbb{F}_{q^k}). $$

By the classical analogies between function fields and number fields, it may be natural to expect that there would be such a cohomological description for the Dedekind  zeta function of a number field $k$, defined by the Euler product
$$ \zeta_k(s) = \prod_{\frak{p} \in {\rm Max}({\cal O}_k)} (1 - {\rm N}(\frak{p})^{-s})^{-1}, $$
over  the set ${\rm Max}({\cal O}_k)$ of maximal ideals of ${\cal O}_k$, where ${\rm N}(\frak{p}) = \#({\cal O}_k/\frak{p})$, and the infinite product converges absolutely for ${\rm Re}(s) > 1$. To be more precise, recall that a proper smooth curve $X$ over a finite field corresponds, by the analogy, to a complete arithmetic curve $\overline{{\rm Spec}({\cal O}_k)} = {\rm Spec}({\cal O}_k) \sqcup \{ \frak{p} | \infty \}$. So we should consider the complete zeta function $\hat{\zeta}_k(s)$ with $\Gamma$-factors.
Due to the analytic nature of $\hat{\zeta}_k(s)$, the would-be  cohomology spaces associated to $\overline{{\rm Spec}({\cal O}_k)}$ must be infinite dimensional. For a smooth proper curve $X$ over $\mathbb{F}_q$, Deninger ([D2]) indeed constructed the ``dynamical" cohomologies $H^i_{\footnotesize \mbox{dyn}}(X)$, which are infinite dimensional $\mathbb{C}$-vector spaces equipped with the action of the operator $\Theta$ induced by an $\mathbb{R}$-action, and he showed the regularized determinant formula for $\zeta_X(s)$:
$$ \zeta_X(s) = \prod_{n=0}^2 {\rm det}_{\infty} (s\;{\rm id} - \Theta \; | \; H^n_{\footnotesize \mbox{dyn}}(X))^{(-1)^{n+1}}, \leqno{(0.2)}$$ 
where ${\rm det}_{\infty}$ stands for the regularized determinant. We note that the $\mathbb{R}$-action (``flow") is analogous to the Frobenius action and closed orbits correspond to Frobenius orbits, namely, closed points of $X$. By the formula (0.2), Deninger suggested that one would attach to a number field $k$ a 3-dimensional space $``\overline{{\rm Spec}({\cal O}_k)}"$ with an $\mathbb{R}$-action and the infinite dimensional cohomology $H^n_{\footnotesize \mbox{dyn}}(``\overline{{\rm Spec}({\cal O}_k)}")$ equipped with the action of the operator $\Theta$ 
 attached to the $\mathbb{R}$-action and one would have the regularized determinant formula
$$ \hat{\zeta}_k(s) = \prod_{n=0}^2 {\rm det}_{\infty} (s\, {\rm id} - \Theta \; | \; H^n_{\footnotesize \mbox{dyn}}(``\overline{{\rm Spec}({\cal O}_k)}"))^{(-1)^{n+1}}. \leqno{(0.3)}$$
We refer to [D13] for the work toward the construction of $``\overline{{\rm Spec}({\cal O}_k)}"$. We also note that there is a different approach to this issue, due to Connes, from noncommutative geometry, where the phase space $``\overline{{\rm Spec}({\cal O}_k)}"$ is defined by a non-commutative quotient of the adeles ([C], [CC]).

Searching for such conjectural cohomology theory, Deninger pointed out the intimate analogies between arithmetic schemes and Riemannian foliated dynamical systems in the series of papers [D6] $\sim$ [D12]. Here a foliated dynamical system is an odd dimensional manifold with one-codimensional foliation and transverse flow, which satisfies some conditions. In particular, a 3-dimensional foliated  dynamical system may be regarded as an analogue of a complete arithmetic curve $\overline{{\rm Spec}({\cal O}_k)}$. Here the closed orbits of the dynamical system,
which are knots in a 3-manifold, correspond to maximal ideals of ${\cal O}_k$. So we note that a 3-dimensional foliated dynamical system picture refines the analogies in arithmetic topology ([M]). As for the conjectural cohomology theory for arithmetic schemes, Deninger noticed that the reduced leafwise cohomologies, which are infinite dimensional in general, should play a similar role for Riemannian foliated dynamical systems, and he gave an important evidence in the theme of explicit formula in analytic number theory.

Now, the purpose of this paper is to give an evidence for Deninger's program by showing an analogue of the regularized determinant formula (0.2)
and the conjectural (0.3), for certain
3-dimensional Riemannian foliated dynamical systems. Let $\frak{S} = (M, g, {\cal F}, \phi)$ be a 3-dimensional Riemannian foliated dynamical system. Namely, $M$ is a $3$-dimensional oriented, connected, closed manifold with Riemannian metric $g$,
${\cal F}$ is a one-codimensional Riemannian foliation on $M$ and $\phi$ is a smooth flow on $M$ such that $\phi$ is transverse to leaves and sends leaves to leaves. Here, the condition on ${\cal F}$ to be Riemannian means that $g$ is bundle-like. In [\'{A}KL1] and [KMNT], such foliated dynamical systems are classified into two types:\\
(i) $M$ is a fiber bundle over $S^1$ whose typical fiber is an oriented, connected, closed surface $S$ and ${\cal F}$ is the bundle foliation, and $\phi$ is leafwise homotopic to the suspension flow of the monodromy map $f : S \rightarrow S$.\\
(ii) ${\cal F}$ is a minimal $\mathbb{R}$-Lie foliation so that all leaves are dense in $M$.

Let ${\cal P}$ be the set of closed orbits of $\phi$. For $\gamma \in {\cal P}$, $\ell(\gamma)$ denotes the length of $\gamma$. A closed orbit $\gamma$ is called simple if ${\rm det}({\rm id} - T_x(\phi^{k\ell(\gamma)}) \; | \; T_x({\cal F})) \neq 0$ for any $x \in \gamma$ and any non-zero integer $k$, where $T_x({\cal F})$ is the tangent space to the foliation at $x$. We assume the condition \\
(A1) any closed orbit is simple, \\
which implies that ${\cal P}$ is a countable set. For a simple closed orbit $\gamma$ and a non-zero integer $k$, we set 
$$\varepsilon_{\gamma}(k) := {\rm sgn} \det ({\rm id} - T_x(\phi^{k\ell(\gamma)}) \; | \; T_x({\cal F})), \;\; \varepsilon_{\gamma} := \varepsilon_{\gamma}(1).$$
Then the zeta function of $\frak{S}$ is defined by the Euler product
$$ \zeta_{\frak{S}}(s) := \prod_{\gamma \in {\cal P}} (1 - e^{-s\ell(\gamma)})^{\epsilon_{\gamma}}$$
over ${\cal P}$. For $x > 0$, let 
$$\nu(x) := \# \{ \gamma \in {\cal P} \; | \; \ell(\gamma) \leq x \}$$
and we assume the condition\\
(A2) $\displaystyle{ a := \limsup_{x \rightarrow \infty} \frac{ \log \nu(x)}{x} }$ exists.\\
It implies that $\zeta_{\frak{S}}(s)$ converges absolutely for ${\rm Re}(s) > a$. 

Let $\bar{H}_{\cal F}^n(M)$ $(0 \leq n \leq 2)$ be the reduced leafwise cohomology on which the infinitesimal operator $\Theta$ induced by $\phi^t$ acts, and we set $\bar{H}^n_{\cal F}(M)_{\mathbb{C}} := \bar{H}^n_{\cal F}(M) \otimes \mathbb{C}$. Then we prove the following regularized determinant formula
$$ \zeta_{\frak{S}}(s) = \prod_{n=0}^{2}  {\rm det}_{\infty} (s\, {\rm id} - \Theta \, | \, \bar{H}^n_{\cal F}(M)_{\mathbb{C}})^{(-1)^{n+1}}. \leqno{(0.4)}$$
for each case of type (i) and type (ii).
\\
Suppose $\frak{S}$ is of type (i) such that the monodromy map $f : S \rightarrow S$ is an Anosov diffeomorphism when $S$ is a torus, more generally, an Axion A diffeomorphism when $S$ is of genus $\geq 1$. Then the conditions (A1), (A2) are satisfied. The formula (0.4) is obtained from the case of the discrete dynamical system of the iteration of $f$ and the classical Lefschetz trace formula with finite dimensional singular cohomologies $H^n(S; \mathbb{R})$. It is similar to obtaining (0.3) from (0.1). 
\\
The case (ii) of dense leaves is more interesting. We assume further the condition (A3) that $\varepsilon_{\gamma}(k) = \varepsilon_{\gamma}$ for any $\gamma \in {\cal P}$ and any $k \in \mathbb{N}$ and the condition (A4) that the flow $\phi$ is leafwise homotopic to a flow which is isometric. We give, in the section 5, examples of Riemannian foliated dynamical systems of type (ii) satisfying (A1) $\sim$ (A4). The key ingredient to prove (0.4) is the distributional dynamical Lefschetz trace formula in [\'{A}K2], [\'{A}KL2], which relates the periodic orbits of the flow $\phi$ with the alternating sum of traces of the infinitesimal generators on the reduced leafwise cohomologies $\bar{H}^n_{\cal F}(M)$. Applying the trace formula, we can relate the dynamical spectral $\xi$-functions with the zeta function, where the dynamical spectral $\xi$ functions are  defined by using the spectra of $\Theta$ as dynamical analogues of the $\xi$-function  attached to the non-trivial zeros of the Riemann zeta function. This process is similar to the method of Deninger in the arithmetic case [D2]. Unlike the arithmetic case, we need to estimate the spectra of $\Theta$ on $\bar{H}^1_{\cal F}(M)$. This is done by using Weyl's asymptotic formula. 

The contents of this paper are organized as follows. The section 1 contains the background materials on Riemannian foliated dynamical systems and leafwise cohomologies. In the section 2, we recall the dynamical Lefschetz trace formula proved in [\'{A}K2], [\'{A}KL2]. In the section 3, we introduce the zeta function of a foliated dynamical system. In the section 4, we recall the regularized determinants and basic properties, where we prove the determinant formula for the zeta function for the case (i). In the section 5, we provide examples of 3-dimensional Riemannian foliated dynamical systems of type (ii) satisfying certain assumptions (A1) $\sim$ (A4). In the section 6, we introduce an analogue of the $\xi$-function for the spectrum of $\Theta$ and then we prove the formula (0.5) for the case (ii) under (A1) $\sim$ (A4), based on the argument of [D2] and  the dynamical Lefschetz trace formula. 
\\
\\
{\it Acknowledgement}. This work was started during the third author's stay at University of Santiago de Compostela on January of 2024. The third author is thankful to the first author and  USC for the support and hospitality. We thank Christopher Deninger for his comments on the draft of this paper. The first author is  partially supported by the grants PID2020-114474GB-I00 (AEI/FEDER, UE) and ED431C 2019/10 (Xunta de Galicia, FEDER). The third author is partially supported by the Japan Society for the Promotion of Science, KAKENHI Grant Number (C) JP22K03270. \\

\begin{center}
{\bf 1. Riemannian foliated dynamical systems and their leafwise cohomologies}
\end{center}

In this section, we provide some notions and properties concerning Riemannian foliated dynamical systems and their leafwise cohomologies, which will be used in the sequel. References are [D6]$\sim$[D12], [Ko], [\'{A}K1], [\'{A}KL1] and [\'{A}KL2].   \\

Let $M$ be an oriented, connected, closed, smooth manifold of dimension $D \geq 3$. A  {\it foliation} ${\cal F}$ of codimension $q$ on $M$ is defined by a family of injectively immersed  $q$-codimensional manifolds $\{ L_{\lambda} \}_{\lambda \in \frak{L}}$  satisfying the following conditions: $M = \bigsqcup_{\lambda \in \frak{L}} L_{\lambda}$ (disjoint union) and $M$ admits a foliated 
atlas $\{ U_{\alpha}, \varphi_{\alpha} \}_{\alpha \in A}$ of codimension $q$, namely, $\varphi_{\alpha}$ is a diffeomorphism from $U_{\alpha}$ to an open ball $B$ in $\mathbb{R}^D$ such that $\varphi_{\alpha}(U_{\alpha} \cap L_{\lambda}) = B \cap (\mathbb{R}^{D - q} \times \{ y_{\alpha} \})$, where $y_{\alpha} \in \mathbb{R}^q$. So a foliated coordinate around $p \in U_{\alpha}$ is written as 
$$\varphi_{\alpha}(p) = (x_{\alpha}, y_{\alpha}), \;\; x = (x_{\alpha}^1, \dots , x_{\alpha}^{D-q}), \; y = (y_{\alpha}^1, \dots , y_{\alpha}^q).$$ 
If $U_{\alpha} \cap U_{\beta} \neq \emptyset$, the change of coordinates should be given in the form  
$$\varphi_{\alpha} \circ \varphi_{\beta}^{-1}(x_{\beta}, y_{\beta})  = (f_{\alpha \beta}(x_{\beta}, y_{\beta}), g_{\alpha \beta}(y_{\beta})),$$ 
where the diffeomorphisms $g_{\alpha \beta} : y_{\beta}(U_{\alpha} \cap U_{\beta}) \rightarrow y_{\alpha}(U_{\alpha} \cap U_{\beta})$ are called  {\it elementary holonomy transformations} of the foliated atlas and form its holonomy cocycle. Each $L_{\lambda}$ is called a {\it leaf} of the foliation ${\cal F}$. A pair $(M, {\cal F})$ is called a {\it foliated manifold}. If every $y_{\alpha}(U_{\alpha})$ is an open subset of a Lie group $G$ and the elementary holonomy transformations $g_{\alpha \beta}$ are restrictions of left translations in $G$, we call the foliation ${\cal F}$ a {\it Lie foliation} with structure group $G$.

A {\it flow} on $M$ is a smooth map $\phi : \mathbb{R} \times M \rightarrow M$ such that $\phi^t := \phi(t, \cdot)$ is a diffeomorphism of $M$ for each $t \in \mathbb{R}$ which satisfies $\phi^0 = {\rm id}_M, \phi^{t+t'} = \phi^t \circ \phi^{t'}$ for any $t, t' \in \mathbb{R}$. A flow $\phi$ on a foliated manifold $(M,{\cal F})$ is called an (${\cal F}$-){\it foliated flow} if $\phi^t$ maps leaves to leaves for any $t \in \mathbb{R}$. We say that foliated flows $\phi$ and $\psi$ are {\it leafwise homotopic} if there is a smooth family of foliated flow $H_s = \{ H_s^t \}$ over $s \in [0,1]$, called a {\it leafwise homotopy}, such that  $H_0^t = \phi^t, H_1^t = \psi^t$ for any $t \in \mathbb{R}$. \\
\\
{\bf Definition 1.1.} We call a triple $(M, {\cal F}, \phi)$ a {\rm foliated dynamical system} if the following conditions are satisfied.\\
(1) $M$ is an oriented, connected, closed, smooth  manifold of odd dimension $D = 2d+1$,\\
(2) ${\cal F}$ is a one-codimensional foliation on $M$ (hence $\dim({\cal F}) = 2d$), \\
(3) $\phi$ is an ${\cal F}$-foliated flow on $M$, and \\
(4) any orbit of $\phi$ is transverse to the leaves.\\
\\
{\bf Remark 1.2.} (1) In [D6]$\sim$[D12] and [Ko], the notion of a foliated dynamical system is defined as a more general object. Namely, instead of the condition (4) in Definition 1.1, we require the weaker condition: There is a finite number of leaves $L_1^{\infty}, \dots , L_r^{\infty}$ such that $\phi^t(L_i^{\infty}) = L_i^{\infty}$ for any $i$ and $t$, and any orbit of $\phi$ is transverse to leaves in $M \setminus \bigcup_{i=1}^r L_i^{\infty}$. The condition (4) is the case that the set of leaves $L_1^{\infty}, \dots , L_r^{\infty}$ preserved by flow is empty. In particular, (4) implies that $\phi$ has no fixed point.  \\
(2) By Definition 1.1, the foliation of a foliated dynamical system $(M, {\cal F}, \phi)$ is a Lie foliation with structure group $\mathbb{R}$.\\
(3) Also, in [D6]$\sim$[D12], the above definition is given for arbitrary dimension $D$. We consider only odd $D=2d+1$ because that is the case where our formulas make sense.\\

In [\'{A}KL1], foliated dynamical systems have been classified as follows (cf. [KMNT] for 3-dimensional case).\\
\\
{\bf Theorem 1.3} ([\'{A}KL1; Theorem 1.1]). {\it Let $(M, {\cal F}, \phi)$ be a foliated dynamical system in the sense of Definition 1.1. Then the foliation ${\cal F}$ is classed as one of the following two types}:\\
(i) {\it  $M$ is a fiber bundle over the circle $S^1$ and ${\cal F}$ is a bundle foliation.} {\it The flow $\phi$ is flow leafwise homotopic to the suspension flow of the monodromy.} \\
(ii) {\it ${\cal F}$ is a minimal $\mathbb{R}$-Lie foliation so that any leaf is dense.}\\
According to the classification in Theorem 1.3, we call a foliated dynamical system $(M, {\cal F}, \phi)$ {\it of type $(i)$} or {\it of type $(ii)$}. \\

Let $(M,{\cal F}, \phi)$ be a foliated dynamical system. Let $T({\cal F})$ be the {\it tangent bundle to the foliation}, which is defined as the sub vector bundle of  $T(M)$, whose total space is the union of the tangent spaces to leaves. Let $X_{\phi}$ be the vector field generated by the flow $\phi$. Then there exists a unique  closed 1-form $\omega_{\phi}$ on $M$, which satisfies 
$$ \omega_{\phi}|_{T({\cal F})} = 0 \;\; \mbox{and}\;\; \omega_{\phi}(X_{\phi}) = 1.$$
Let $T_0(M)$ be the rank one vector bundle over $M$ generated by $X_{\phi}$. We thus have the decomposition $ T(M) = T({\cal F}) \oplus T_0(M).$

Let $(M, {\cal F}, \phi)$ be a foliated manifold and suppose that $M$ is a Riemannian manifold with metric $g$. Then the foliation ${\cal F}$ is called a {\it Riemannian foliation} if the metric $g$ is {\it bundle-like}, namely,  $X_{\phi}$ is of constant norm (say of norm one) and orthogonal to $T({\cal F})$, $T(M) = T({\cal F}) \bot T_0(M)$. \\
\\
{\bf Definition 1.4.} We call a 4-tuple $(M, g, {\cal F}, \phi)$ a {\rm Riemannian foliated dynamical system} if the following conditions are satisfied.\\
(1) $M$ is an oriented, connected, closed, Riemannian manifold of odd dimension $2d+1$ with metric $g$,\\
(2) ${\cal F}$ is a one-codimensional Riemannian foliation on $M$ (hence $\dim({\cal F}) = 2d$), \\
(3) $\phi$ is an ${\cal F}$-foliated flow on $M$, and \\
(4) any orbit of $\phi$ is transverse to leaves.\\
\\
Though the following lemma might be well known, we give a proof for the sake of readers.\\
\\
{\bf Lemma 1.5.} {\it Let $(M, {\cal F}, \phi)$ be a foliated dynamical system. Then there is a Riemannian metric $g$ on $M$ such that $(M, g, {\cal F}, \phi)$ is a Riemannian foliated dynamical system. }\\
\\
{\it Proof.} Consider the Euclidean structure $g_0$ on $T_0(M)$ determined by requiring that $X_\phi$ is of norm one at every point. Let $g_{\mathcal F}$ be any Euclidean structure on $T({\mathcal F})$. Then there is a unique Riemannian metric $g$ on $M$ whose restrictions to $T_0(M)$ and $T({\mathcal F})$ are $g_0$ and $g_{\mathcal F}$, respectively,  
and such that $T_0(M)$ and $T({\mathcal F})$ are orthogonal. Such a metric $g$ is bundle-like for $\mathcal F$. $\;\; \Box$
\\

Now let $(M, {\cal F}, \phi)$ be a  foliated dynamical system. Leafwise differential forms of degree $n$  are defined as the smooth sections of $\wedge^n T^*({\cal F})$,
$$ {\cal A}_{\cal F}^n(M) := \Gamma(M, \wedge^n T^*({\cal F})),$$
and the exterior derivatives along the leaves 
$$ d_{\cal F}^n : {\cal A}_{\cal F}^n(M) \longrightarrow {\cal A}_{\cal F}^{n+1}(M)$$
satisfy $d_{\cal F}^{n+1} \circ d_{\cal F}^n = 0$. Thus we have the de Rham complex along the leaves $ ({\cal A}_{\cal F}^{\bullet}, d_{\cal F}^{\bullet})$
and  we can form the {\it leafwise cohomology} of ${\cal F}$:
$$ H^n_{\cal F}(M) := {\rm Ker}(d_{\cal F}^n)/{\rm Im}(d_{\cal F}^{n-1}).$$
For our purposes, these cohomology groups are too subtle (they may not be Hausdorff) and so we consider the {\it reduced leafwise cohomology}
$$ \bar{H}_{\cal F}^n(M) := {\rm Ker}(d_{\cal F}^n)/\overline{{\rm Im}(d_{\cal F}^{n-1})}, $$
where  the quotient is taken by the topological closure of ${\rm Im}(d_{\cal F}^{n-1})$ in the natural Fr\'{e}chet topology on ${\cal A}^n_{\cal F}(M)$. 
Note that reduced leafwise cohomology groups may be infinite dimensional, even if the leaves are dense ([AH]). The exterior product of the leafwise differential forms induces the cup product pairing
$$ \cup :  \bar{H}_{\cal F}^m(M) \times \bar{H}_{\cal F}^n(M) \longrightarrow \bar{H}_{\cal F}^{m+n}(M),$$
which makes $\bar{H}_{\cal F}^{\bullet}(M)$ into a graded $\bar{H}_{\cal F}^0(M)$-algebra.

For a smooth map $f : M \rightarrow N$ between foliated manifolds $(M, {\cal F}_M)$ and $(N, {\cal F}_N)$, which sends leaves in ${\cal F}_M$ into leaves in ${\cal F}_N$, continuous pullback maps
$$ f^* : {\cal A}_{{\cal F}_N}^n(N) \longrightarrow {\cal A}_{{\cal F}_M}^n(M)$$
are defined for any $n$, and they commute with $d_{\cal F}$ and respect the exterior product of differential forms along the leaves. Hence we have a continuous homomorphism
of reduced leafwise cohomology algebras
$$ f^{\bullet} : \bar{H}_{{\cal F}_N}^{\bullet}(M) \longrightarrow \bar{H}_{{\cal F}_M}^{\bullet}(N).$$

Since $\phi^t : M \rightarrow M$ sends leaves into leaves, $\phi^{t*}$ induces an $\mathbb{R}$-linear automorphism $\phi^{t*}$ of $\bar{H}_{\cal F}^{n}(M)$ for each $n$.
If foliated flows $\phi$ and $\psi$ on $(M,{\cal F})$ are leafwise homotopic, then we have $\phi^{t*} = \psi^{t*}$ as $\mathbb{R}$-linear automorphisms of $\bar{H}_{\cal F}^{n}(M)$ for any $n$  ([\'{A}KL2; 3.2.16]). Let 
$$\Theta : \bar{H}_{\cal F}^n(M) \longrightarrow \bar{H}_{\cal F}^n(M)$$
denote the infinitesimal generator of $\phi^{t*}$ defined by
$$\Theta(h) = \lim_{t \rightarrow 0} \frac{ \phi^{t*}(h) - h}{t} \;\;\;\; \mbox{for} \; h \in \bar{H}_{\cal F}^n(M)$$ 
\\
{\bf Example 1.6.} Let $(M, {\cal F}, \phi)$ be a 3-dimensional foliated dynamical system of type (i) so that $M$ is a fiber bundle over $S^1$, $\pi : M \to S^1$, whose typical fiber is  an oriented, connected, closed surface $S$, ${\cal F}$ is the bundle foliation whose leaves are fibers $\pi^{-1}(s) \; (s \in S^1)$, and $\phi$ is leafwise homotopic to the suspension flow of a diffeomorphism $f: S \rightarrow S$. The surface $S$ admits a Riemannian metric, which induces a bundle-like Riemannian metric $g$ on $M$ such that ${\cal F}$ is a Riemannian foliation. (Actually, ${\cal F}$ admits a structure of complex foliation.) Then $\bar{H}^n_{\mathcal F}(M)$ is isomorphic to the smooth sections of a flat vector bundle over $S^1$ whose typical fiber is $H^n(S,{\mathbb R})$. This flat vector bundle is the suspension of $H^n(f)$. The derivative of the smooth sections of this flat vector bundle with respect to the suspension flow defines the infinitesimal operator $\Theta$ on $\bar{H}^n_{\cal F}(M)$. 
\\

Next we suppose that $(M, g, {\cal F}, \phi)$ is a Riemannian foliated dynamical system. Let $\langle \cdot , \cdot \rangle_{\cal F}$ be the Riemannian metric on $\wedge^{\bullet} T^*({\cal F})$ induced by $g$. Then we have the inner product on ${\cal A}_{\cal F}^{\bullet}(M)$ defined by
$$ (\alpha, \beta)_{\cal F} := \int_M \langle \alpha, \beta \rangle_{\cal F} {\rm vol} \;\;\;\; (\alpha, \beta \in {\cal A}_{\cal F}^{\bullet}), $$
 where ${\rm vol}$ is the volume form on $M$ induced by $g$.  Let
 $$ \Delta_{\cal F} = d_{\cal F} \circ d_{\cal F}^{*} + d_{\cal F}^{*} \circ d_{\cal F}$$
 denote the leafwise Laplacian, where $d_{\cal F}^*$ is the formal adjoint of $d_{\cal F}$ on $M$. Since the metric $g$ is bundle-like, the formal adjoint on the leaves is the same as the formal adjoint on the ambient manifold. We assume that $T({\cal F})$ is orientable. Choosing an orientation,  via $\langle \cdot , \cdot \rangle_{\cal F}$, we have a volume form in ${\cal A}_{\cal F}^{2d}(M)$ and hence a Hodge $*$-operator
 $$ *_{\cal F} :  \wedge^n T_x^*({\cal F}) \stackrel{\sim}{\longrightarrow} \wedge_x^{2d-n} T_x^*({\cal F}) \;\;\;\; \mbox{for any}\; x \in M,$$
which is determined by the condition that
$$ v \wedge *_{\cal F} w = \langle v, w \rangle_{\cal F} {\rm vol}_{{\cal F},x} \;\;\;\; \mbox{for any}\; v, w \in \wedge^{\bullet} T_x^*({\cal F}).$$
These fiberwise $*$-operators induce the leafwise $*$-operator on differential forms along leaves:
$$ *_{\cal F} : {\cal A}_{\cal F}^n(M) \stackrel{\sim}{\longrightarrow} {\cal A}_{\cal F}^{2d-n}(M).$$
The following result proved in [AK1] may be regarded as a leafwise analogue of the Hodge theorem for leawise cohomology.\\
\\
{\bf Theorem 1.7} ([\'{A}K1; Theorem B, Corollary C]). {\it The natural map }:
$$  {\rm Ker}(\Delta_{\cal F}^n) \stackrel{\sim}{\longrightarrow} \bar{H}_{\cal F}^n(M); \;\; \omega \mapsto \omega \; \mbox{mod}\;  \overline{{\rm Im}(d_{\cal F})}.$$
{\it is a topological isomorphism of Fr\'{e}chet spaces.} We denote its inverse by ${\cal H} : \bar{H}_{\cal F}^n({\cal F}) \rightarrow {\rm Ker}(\Delta_{\cal F}^n)$.  \\

Since the leafwise Hodge $*$-operator commutes with $\Delta_{\cal F}$, it induces the isomorphism
$$ *_{\cal F} : {\rm Ker}(\Delta_{\cal F}^n) \stackrel{\sim}{\longrightarrow} {\rm Ker}(\Delta_{\cal F}^{2d-n}).$$
By the above leafwise Hodge theorem, we therefore obtain isomorphisms
$$ *_{\cal F} : \bar{H}_{\cal F}^n(M) \stackrel{\sim}{\longrightarrow} \bar{H}_{\cal F}^{2d-n}(M) \;\;\;\; \mbox{for all}\; n$$
We define the {\it trace map} 
$$ {\rm tr} : \bar{H}_{\cal F}^{2d}(M) \longrightarrow \mathbb{R}$$
by the formula
$$ {\rm tr}(h) = \int_M *_{\cal F}(h) {\rm vol} := \int_{M} *_{\cal F}({\cal H}(h)) {\rm vol}.$$
By the leafwise Hodge theorem above, the inner product $(\cdot , \cdot)_{\cal F}$ on ${\cal A}_{\cal F}^{\bullet}(M)$ induces the inner product on the reduced leafwise cohomology
$\bar{H}_{\cal F}^{\bullet}(M)$, which is described by the trace map:
$$ \begin{array}{ll} (h, h')_{\cal F} & = \displaystyle{\int_M \langle {\cal H}(h), {\cal H}(h') \rangle_{\cal F} {\rm vol} } \\
                                      & = {\rm tr}(h \cup {\cal *}_{\cal F} h'). 
\end{array}$$

We say that the flow $\phi$ is {\it isometric} on the leaves of ${\cal F}$ if for any $x \in M$, $v, w \in T_x({\cal F})$ and $t \in \mathbb{R}$, the following equality holds:
$$ g(T_x(\phi^t)(v), T_x(\phi^t)(w)) = g(v,w).$$
On leafwise differential forms, it is equivalent to say that for any $\alpha, \beta \in {\cal A}_{\cal F}^{\bullet}$, the following equality holds:
$$ (\phi^{t*}(\alpha), \phi^{t*}(\beta))_{\cal F} = (\alpha,\beta)_{\cal F}.$$
\\
{\bf Theorem 1.8.} {\it Let $(M, g, {\cal F}, \phi)$ be a $3$-dimensional Riemannian foliated dymanical system of type $(ii)$ $(\mbox{cf. Theorem 1.3})$ such that $\phi$ is leafwise homotopic to a foliated flow $\psi$ which is isometric on the leaves. Let $\Theta$ be the infinitesimal generator of $\phi^{t*}$. Then we have}
$$ \bar{H}_{\cal F}^0(M)
\simeq \bar{H}_{\cal F}^2(M) \simeq \mathbb{R} \;\; \mbox{and}\;\; \Theta = 0,$$
{\it and}
$$ \Theta \; \mbox{{\it is skew-symmetric on}}\; \bar{H}_{\cal F}^1(M) \; \mbox{{\it with respect to}} \; (\cdot, \cdot)_{\cal F}.$$
\\
{\it Proof.} Since the leafwise Hodge $*_{\cal F}$-operator commutes with $\psi^{t*}$, the leafwise Hodge Theorem 1.7 yields $(\bar{H}_{\cal F}^0(M), \Theta) \simeq (\bar{H}_{\cal F}^2(M), \Theta)$. Since $\phi^{t*} = \psi^{t*}$, the remaining assertions follow from [D8; Theorem 2.1]. $\Box$\\

\begin{center}
{\bf 2. Distributional dynamical Lefschetz trace formulas}
\end{center}

In this section, we recall the distributional dynamical Lefschetz trace formulas for Riemannian foliated dynamical systems proved in [AK2]. The formula relates the periodic orbits of the flow with the alternating sum of traces of the infinitesimal generators on the leafwise cohomologies.  \\

Let $(M, {\cal F}, \phi)$ be a foliated dynamical system. Let ${\cal P}$ denote the set of all closed orbits of the flow $\phi$. For $\gamma \in {\cal P}$,
let $\ell(\gamma)$ denote the {\it length} of $\gamma$, the minimal $t>0$ such that $\phi^t(p) = p$ for some $p \in \gamma$. It is independent of the choice of $p \in \gamma$. So the length $\ell(\gamma)$ is also defined by the isomorphism $\mathbb{R}/\ell(\gamma)\mathbb{Z} \stackrel{\sim}{\rightarrow} \gamma; \; t \mapsto \phi^t(x)$,    for any $x\in\gamma$. 
A closed orbit  $\gamma \in {\cal P}$ is called {\it simple} or {\it non-degenerate} if the following condition holds: For any $x \in \gamma$ and any non-zero integer $k$   the automorphism $T_x(\phi^{k\ell(\gamma)})$ of $T_x(M)$ has the eigenvalue 1 with algebraic multiplicity one, equivalently, $\det ({\rm id} - T_x(\phi^{k\ell(\gamma)}) \; | \; T_x({\cal F})) \neq 0$. We say a flow $\phi$ is {\it simple} if any closed orbit is simple.
 For a simple closed orbit $\gamma$ and non-zero integer $k$, we set
$$\varepsilon_{\gamma}(k) := {\rm sgn} \det ({\rm id} - T_x(\phi^{k\ell(\gamma)}) \; | \; T_x({\cal F})), \;\; \varepsilon_{\gamma} := \varepsilon_{\gamma}(1). \leqno{(2.1)}$$
We call $\varepsilon_{\gamma}(k)$ (resp. $\varepsilon_{\gamma}$) the {\it index} of $(\gamma, k)$ (resp. of $\gamma$). \\
\\
{\bf Lemma 2.2.} {\it Let $(M, {\cal F}, \phi)$ be a foliated dynamical system such that the flow $\phi$ is simple. Then the set of closed orbits ${\cal P}$ is a countable subset of $\mathbb{R}_+ = (0,\infty)$; in fact, there are finitely many closed orbits with periods in any compact interval of $\mathbb R_+$.} \\
\\
{\it Proof.}  Since any closed orbit is simple, we see that the following two maps
$$ M \times \mathbb{R}_+ \longrightarrow M \times M \times \mathbb{R}_+; \;\; (p, t) \mapsto (p, \phi^t(p), t)$$
and
$$ M \times \mathbb{R}_+ \longrightarrow M \times M \times \mathbb{R}_+ ; \;\; (p, t) \mapsto (p, p, t)$$
are transversal. Then the closed orbits are isolated if we take the periods in a compact interval $I = [a, b]$ of $\mathbb{R}_+$. Therefore there are finitely many closed  orbits with periods in $I$.  $\;\;\Box$\\
\\
{\bf Lemma 2.3.} ${\rm Inf} \{ \ell(\gamma) \; | \; \gamma \in {\cal P} \} > 0.$ \\
\\
{\it Proof}. By the proof of Lemma 2.2, there are finitely many closed orbits with periods in a compact interval $I$ of $\mathbb{R}_+$. This shows that the set $\{ \ell(\gamma) \; | \; \gamma \in {\cal P} \}$ is a discrete subset of $\mathbb{R}_+$. Its infimum can not be $0$, otherwise, there would be some fixed points of the flow by the compactness of $M$. Hence the assertion follows. $\;\; \Box$ \\

Let $(M, g, {\cal F}, \phi)$ be a Riemannian foliated dynamical system. By Theorem 1.7, the leafwise cohomologies $\bar{H}_{\cal F}^n(M)$ are equipped with the inner product $(\cdot , \cdot)_{\cal F}$ so that $(\bar{H}_{\cal F}^n(M),(\cdot , \cdot)_{\cal F})$ are pre-Hilbert spaces. Let $\hat{H}_{\cal F}^n(M)$ denote their Hilbert space completions. For an open set $U$ of $\mathbb{R}$, let ${\cal D}(U) = C_{\rm c}^{\infty}(U)$ denote the space of compactly supported smooth functions on $U$ and ${\cal D}'(U)$ be the 
space of distributions on $U$. For $a \in U$, let $\delta_a \in {\cal D}'(U)$ denote the Dirac distribution at $a$. Then we have the following distributional dynamical Lefschetz trace formula. \\
\\
{\bf Theorem 2.4.} {\it Let $(M, g, {\cal F}, \phi)$ be a Riemannian foliated dynamical system.}\\
(1) ([\'{A}K2; Theorem 1.1]). {\it For any $t \in \mathbb{R}$, the linear operator $\phi^{t*}$ is bounded on $\hat{H}_{\cal F}^n(M)$ for each $n$. For every $f \in {\cal D}(\mathbb{R})$, the operator} 
$$ A_f := \int_{\mathbb{R}} f(t) \phi^{t*} dt$$
{\it on $\hat{H}_{\cal F}^n(M)$ is of trace class, and the functional $f \mapsto {\rm tr}(A_f)$ defines a distribution on $\mathbb{R}$. We set}:
$$ {\rm tr}(\phi^* \; | \; \bar{H}_{\cal F}^n(M))(f) = {\rm tr}(A_f).$$
(2) ([\'{A}K2; Theorem 1.3]). {\it Assume that the flow $\phi$ is simple. Then the following formula holds in ${\cal D}'(\mathbb{R_+})$}:
$$ \displaystyle{   \sum_{n=0}^{2d} (-1)^n  {\rm tr}(\phi^* \; | \; \bar{H}_{\cal F}^n(M)) =  \sum_{\gamma \in {\cal P}} \ell(\gamma) \sum_{k=1}^{\infty} \varepsilon_{\gamma}(k) \delta_{k\ell(\gamma)}. } $$
\vspace{0.2cm}\\
Let $\bar{H}^{n}_{\cal F}(M)_{\mathbb{C}} := \bar{H}^n_{\cal F}(M) \otimes_{\mathbb{R}} \mathbb{C}$ and $\hat{H}^n_{\cal F}(M)_{\mathbb{C}} := \hat{H}_{\cal F}^n(M) \otimes_{\mathbb{R}} \mathbb{C}$, which are equipped with Hermitian extension of $(\cdot , \cdot)_{\cal F}$.  By Theorems 1.8 and 2.4, we have the following. \\
\\
{\bf Theorem 2.5.} {\it Let $(M,g, {\cal F}, \phi)$ be a $3$-dimensional Riemannian foliated dynamical system of type $(ii)$ $(\mbox{cf. Theorem 1.3})$. Assume that $\phi$ is simple and that $\phi$ is leafwise homotopic to a foliated flow which is isometric.} \\
(1)  {\it $\Theta$ has pure point spectrum ${\rm Sp}^1(\Theta)$ in $\hat{H}^1_{\cal F}(M)_{\mathbb{C}}$ which is discrete in $\mathbb{R} \sqrt{-1}$, and $\hat{H}^1_{\cal F}(M)_{\mathbb{C}}$ is a direct sum of $\Theta$-invariant subspaces and each $\alpha \in {\rm Sp}^1(\Theta)$ occurs with finite multiplicities.} \\
(2) {\it According to (1), the elements of $\operatorname{Sp}^1(\Theta)$ form a sequence $\rho_0,\rho_1.\dots$, taking multiplicities into account, with $|\rho_j| \le |\rho_{j+1}|$ for all $j$. Then there is some $C>0$ such that $|\rho_j| \ge Cj^{1/3}$.}\\
(3) {\it We have the following equality in $\mathcal D'(\mathbb R_+)$}:
$$ 2 - \sum_{\rho \in {\rm Sp}^1(\Theta)} e^{\rho x} = \sum_{\gamma \in {\cal P}} \ell(\gamma) \sum_{k = 1}^{\infty} \varepsilon_{\gamma}(k)  \delta_{k \ell(\gamma)}. $$
\\
{\it Proof.} (1) Since $\Theta$ is skew-symmetric (Theorem 1.8), it has purely imaginary eigenvalues. 
We can assume  $\phi$ is isometric on the leaves because leafwise homotopic foliated flows induce the same action on the reduced leafwise cohomology. Then  $\Delta_M$ preserves $\operatorname{Ker} \Delta_{\mathcal F}$ in ${\mathcal A}_{\mathcal F}^1(M)$, where we have $\Delta_M = -\Theta^2$. So the spectrum of the operator $-\Theta^2$ in $\hat{H}^1_{\mathcal F}(M)_{\mathbb C}$ is discrete and contained in $[0,\infty)$. Hence the spectrum of $\Theta$ in $\hat{H}^1_{\mathcal F}(M)_{\mathbb C}$ is also discrete; from here we can also obtain that it is contained in ${\mathbb R}\sqrt{-1}$. \\
(2) Let $0\le \lambda_0 \le \lambda_1 \le \cdots$ denote the spectrum of $\Delta_M$ on $\mathcal A^1(M)$, taking multiplicities into account. Since $\dim M = 3$, the Weyl's asymptotic formula states that $\lambda_k \sim C_0k^{2/3}$ for some $C_0>0$ (for example, see Corollary 2.43 of [BGV]). Then the property stated in (2) follows because $-\rho_j^2$ is the subsequence of $\lambda_k$ given by the eigenvalues whose eigenvectors belong to $\operatorname{Ker}\Delta_{\mathcal F}$ in $\mathcal A^1_{\mathcal F}(M)$.\\
(3) This equality follows from (1) and Theorem 2.4.  $\;\;\Box$\\
\\
{\bf Remark 2.6.} By Lemmas 2.2 and 2.3, the right-hand side of the equalities of (2) in Theorem 2.4 and (3) in Theorem 2.5 is a Radon measure on $\mathbb{R}_+$ supported in $[m,\infty)$, where $m>0$ is the infimum stated in Lemma 2.3. Its left-hand side is a series of Radon measures on $\mathbb R_+$. However this series might not be convergent in the space of Radon measures on $\mathbb R_+$, even though it is convergent in the space of ditributions on $\mathbb R_+$. 
\\

\begin{center}
{\bf 3. Dynamical zeta functions}
\end{center}

Let $\frak{S} = (M, {\cal F}, \phi^t)$ be a  foliated dynamical system such that $\phi$ is simple. 
By Lemma 2.2, the set ${\cal P}$ of closed orbits of $\phi$ is countable. Recall that for $\gamma \in {\cal P}$, $\ell(\gamma)$ denotes the length of $\gamma$.
We then define the {\it zeta function} of $\frak{S}$ by
$$ \zeta_{\frak{S}}(s) := \prod_{\gamma \in {\cal P}} (1 - e^{-s\ell(\gamma)})^{-\varepsilon_{\gamma}},$$
where $\varepsilon_{\gamma}$ is the index of $\gamma$ defined in (2.1). We also introduce two zeta functions by 
$$ \zeta_{\frak{S}}^{+}(s) = \prod_{\gamma \in {\cal P}} (1 - e^{-s\ell(\gamma)})^{-1}, \;\; \zeta_{\frak{S}}^{-}(s) :=  \prod_{\gamma \in {\cal P}} (1 - e^{-s\ell(\gamma)}),$$
which are the special cases, according as every $\varepsilon_{\gamma}$ is $-1$ or $+1$.\\
\\
{\bf Lemma 3.1.} {\it For $r > 0$, the following conditions are equivalent.}\\
(1) {\it $\zeta_{\frak{S}}^{+}(s)$ converges absolutely for ${\rm Re}(s) > r$.}\\
(2) {\it $\zeta_{\frak{S}}^{-}(s)$ converges absolutely for ${\rm Re}(s) > r$.}\\
(3) {\it $\displaystyle{ \sum_{\gamma \in {\cal P}} e^{-{\rm Re}(s) \ell(\gamma)} }$ converges  for ${\rm Re}(s) > r$.}\\
{\it Let $\nu(x) := \# \{ \gamma \in {\cal P} \; | \; \ell(\gamma) \leq x \}$ for $x > 0$ and consider the following condition}:\\
(4)  {\it $\displaystyle{ a := \limsup_{x \rightarrow \infty} \frac{ \log \nu(x)}{x} }$ exists.}\\
{\it Then $(4)$ implies $(1)$, $(2)$ and $(3)$ with $r = a$.}\\
\\
{\it Proof}. This lemma follows from the proof of [S; Proposition II-2-2].  $\;\; \Box$\\
\\
{\bf Proposition 3.2.} {\it Assume that the equivalent conditions $(1), (2), (3)$ in Lemma 3.1 hold. Then the zeta function $\zeta_{\frak{S}}(s)$ converges absolutely for ${\rm Re}(s) > {\rm max}\{ r, \frac{\log 2}{m} \} $, where $m := {\rm Inf} \{ \ell(\gamma) \; | \; \gamma \in {\cal P} \}$ (cf. Lemma 2.3).}\\
\\
{\it Proof}. Write
$$ (1 - e^{-s\ell(\gamma)})^{-\varepsilon_{\gamma}} = 1 + u_{\gamma}(s), \;\; \gamma \in {\cal P}.$$
By the standard fact on an infinite product, it suffices to show that the series 
$$\displaystyle{ \sum_{\gamma \in {\cal P}} |u_{\gamma}(s)|}$$
converges for ${\rm Re}(s) > {\rm max}\{ r, \frac{\log 2}{m}\}$.\\
When $\varepsilon_{\gamma} = -1$, we have $u_{\gamma}(s) = -e^{-s\ell(\gamma)}$ and so 
$$|u_{\gamma}(s)| = e^{-{\rm Re}(s) \ell(\gamma)} \;\; \mbox{for} \; \gamma \in {\cal P}.$$
When $\varepsilon_{\gamma} = +1$, we have 
$$\displaystyle{u_{\gamma}(s) = \sum_{n=1}^{\infty} e^{- n s \ell(\gamma)}.}$$
If $ {\rm Re}(s) \geq \frac{\log 2}{m}$, we have
$$ |u_{\gamma}(s)| \leq \sum_{n=1}^{\infty} |e ^{-n s \ell(\gamma)}| = \sum_{n=1}^{\infty} e ^{-n {\rm Re}(s) \ell(\gamma)} = \frac{ e^{-{\rm Re}(s) \ell(\gamma)}}{1 -  e^{-{\rm Re}(s) \ell(\gamma)}}\leq 2 e^{-{\rm Re}(s)\ell(\gamma)}$$
for any $\gamma \in {\cal P}$. By the assumption (3),  $\displaystyle{ \sum_{\gamma \in {\cal P}} 2e^{-{\rm Re}(s) \ell(\gamma)} }$ converges  for ${\rm Re}(s) > r$. Hence
 $\displaystyle{ \sum_{\gamma \in {\cal P}} |u_{\gamma}(s)|}$ converges for ${\rm Re}(s) \geq {\rm max}\{r, \frac{\log 2}{m} \}$. $\;\; \Box$ \\
\\
{\bf Example 3.3.} Let $(M, g, {\cal F}, \phi)$ be a 3-dimensional Riemannian foliated dynamical system of type (i) so that $M$ is a fiber bundle over $S^1$ with fiber being an oriented, connected, closed Riemannian surface $S$, ${\cal F}$ is the bundle foliation, and $\phi$ is leafwise homotopic to the suspension flow of the monodromy diffeomorphism $f : S \rightarrow S$. When $S$ is a torus, we assume that $f$ is an Anosov diffeomorphism. When the genus of $S$ $\geq 2$, we assume that $f$ is an Axiom A diffeomorphism. Then, for any periodic point $p\in S$ of $f$ with period $N$, $T_p(S)$ is decomposed into the stable vector subspace $T^{\text{\rm s}}_p(S)$ and the unstable vector subspace $T^{\text{\rm us}}_p(S)$. In the case where both $T^{\text{\rm s}}_p(S)$ and $T^{\text{\rm us}}_p(S)$ are of dimension one, the eigenvalue $\lambda$ of $T_p(f^N)$ on $T^{\text{\rm s}}_p(S)$ is positive and smaller than $1$ and the eigenvalue $\mu$ of $T_p(f^N)$  on $T^{\text{\rm us}}_p(S)$ is bigger than $1$. It follows that, for any closed orbit $\gamma$ of $\phi$ and any point $p\in\gamma$, $T_p({\mathcal F})$ is decomposed into the stable vector subspace $T^{\text{\rm s}}_p({\mathcal F})$ and the unstable vector subspace $T^{\text{\rm us}}_p({\mathcal F})$, and the eigenvalue $\lambda$ of $T_p(\phi^{\ell(\gamma)})$ on $T^{\text{\rm s}}_p({\cal F})$
is positive and smaller than $1$ and the eigenvalue $\mu$ of $T_p(\phi^{\ell(\gamma)})$  on $T^{\text{\rm us}}_p({\mathcal F})$ is bigger than $1$. Therefore $\det({\rm id} - \phi^{k\ell(\gamma)} \; | \; T({\cal F})) = (1 - \lambda^k)(1- \mu^{k}) < 0$ and so $\varepsilon_{\gamma}(k) = -1$ for all $k \neq 0$. In the cases where $T^{\text{\rm s}}_p(S)$ or $T^{\text{\rm us}}_p(S)$ is equal to $T_p(S)$, we similarly get that $\epsilon_\gamma(k) = 1$ for all $k$. It is known that for Anosov diffeomorphisms of a torus or, more generally, an Axiom A diffeomorphism of any surface, the growth of $\nu(x)$ is at most exponential in $x$. In fact, the growth rate of the periodic orbits in this case is given by the expression involving the topological entropy ([Bo]). Hence, by Lemma 3.1, $\zeta_{\mathfrak S}(s)$ converges absolutely for $\operatorname{Re}(s) > r$ with some $r > 0$. For an Anosov diffeomorphism on the torus, only the case $\epsilon_\gamma(k) = -1$ can happen for all $\gamma$ and $k$, obtaining $\zeta_{\mathfrak S}(s) = \zeta^-_{\mathfrak S}(s)$.\\

Let $\mu$ be the Radon measure on $\mathbb R_+$ given by the right-hand side of (2) of Theorem 2.4 and (3) in Theorem 2.5 (see Remark 2.6); i.e.,
$$
\mu = \sum_{\gamma\in\mathcal P} \ell(\gamma) \sum_{k=1}^\infty \delta_{k\ell(\gamma)}.
$$
\\
{\bf Proposition 3.4.}  {\it We assume that the condition (4) of Lemma 3.1 is true and let $m := {\rm Inf} \{ \ell(\gamma) \; | \; \gamma \in {\cal P} \}$.}\\
(1) {\it We have
$$
\#\{(\gamma,k)\in{\mathcal P} \times {\mathbb Z}_+ \mid k\ell(\gamma) \le x\}
\le \frac{xe ^{ax}}{m}+1.
$$
for all $x\gg 0$. }\\
(2) {\it If a continuous function $f: \mathbb{R}_+ \to {\mathbb C}$ satisfies $|f(x)|\le e^{-bx}$ for some $b>a$ and all $x\gg 0$, then $f$ is absolutely integrable with respect to $|\mu|$.}\\
\\
{\it Proof.} (1) We can easily see this inequality by (4) of Lemma 3.1 and the definition of $m$. \\
(2) Take $x_0 \gg 0$ so that the above inequality holds for $x\ge x_0$. Then we have
\begin{align*}
\int_0^\infty |f(x)|\,d|\mu|(x)
&\le \sum_{\gamma\in{\mathcal P}}\ell(\gamma) \sum_{k=1}^\infty  |f(k\ell(\gamma))| \\
&\le \sum_{\gamma\in{\mathcal P},\ \ell(\gamma)\le x_0}
\ell(\gamma) \sum_{k=1}^{\lfloor x_0/\ell(\gamma) \rfloor}  
|f(k\ell(\gamma))| \\
&\phantom{={}}{} + \sum_{x=1}^\infty x \left(\frac{xe^{ax}}{m} + 1\right) |f(x)| \\
&\le \sum_{\gamma\in{\mathcal P},\ \ell(\gamma)\le x_0}
\ell(\gamma) \sum_{k=1}^{\lfloor x_0/\ell(\gamma) \rfloor}  
|f(k\ell(\gamma))| \\
&\phantom{={}}{} + \sum_{x=1}^\infty
\left( \frac{x^2e^{ax}}{m} + x \right) e^{-b(x-1)} \\
&<\infty. \;\; \Box
\end{align*}
\\

\begin{center}
{\bf 4. Regularized products and regularized determinants}
\end{center}

In this section, we recall the definitions and properties of regularized products and regularized determinants, following [D1; $\S$ 1,2], [D4; $\S$ 1,2] and [D11; $\S$ 3].\\

Let $( \alpha_{\nu})$ $(\nu = 1, 2, \dots)$ be a finite or infinite sequence of complex numbers $\alpha_{\nu}$ with given argument $-\pi < {\rm Arg}(\alpha_{\nu}) \leq \pi$, where we assume that only finitely many $\alpha_{\nu}$ are zero. We say that the {\it regularized product} $\underline{\prod} \alpha_{\nu}$ is defined if the following  holds: 

 Let $N \geq 0$ be such that $\alpha_{\nu} \neq 0$ for $\nu > N$. Then the Dirichlet series 
$\displaystyle{\sum_{\nu > N} \alpha_{\nu}^{-s}}$ with $\alpha_{\nu}^{-s} = |\alpha_{\nu}|^{-s} e^{-\sqrt{-1}s{\rm Arg}(\alpha_{\nu})}$ converges absolutely for 
${\rm Re}(s) > r$ with some $r > 0$, and has the analytic continuation to a 
holomorphic function $Z_N(s)$ for ${\rm Re}(s) > - \delta$ with some $\delta > 0$. 

Assuming this, we define the {\it regularized product} of $(\alpha_{\nu})$ by
$$   \underline{\prod}_{\nu} \alpha_{\nu} := \left(   \prod_{\nu = 0}^N \alpha_{\nu} \right) \exp(- Z_N'(0),)$$ \\
which is independent of the choice of $N$. \\

Let $V$ be a complex vector space of countable dimension and let $\theta$ be an endomorphism of $V$. We say that the {\it regularized determinant} ${\rm det}_{\infty}(\theta)$ is defined if the following holds:\\
$\bullet$  $V$ is a direct sum of finite dimensional $\theta$-invariant subspaces and the 

eigenvalues $\alpha$ of $\theta$ occur with finite multiplicities $m(\alpha)$.\\
Assuming this, let ${\rm Sp}(\theta \,|\, V)$ or simply ${\rm Sp}(\theta)$ denote the set of eigenvalues of $\theta$ taken the multiplicities into account. Then we may arrange ${\rm Sp}(\theta)$ in a sequence $(\alpha_{\nu})$ such that each eigenvalue $\alpha$ occurs exactly $m(\alpha)$ times. Given a choice of the argument of $\alpha_{\nu}$ with $-\pi < {\rm Arg}(\alpha_{\nu}) < \pi$ for each $\alpha_{\nu}$, we define the  {\it regularized determinant} of $\theta$ by
$$   {\rm det}_{\infty}(\theta \; | \; V) := \underline{\prod}_{\nu} \alpha_{\nu},$$ 
which is independent of the ordering of the eigenvalues in ${\rm Sp}(\theta)$.\\
\\
The formulas in the following example were shown in [D1], [D2] by using the Hurwitz zeta function. \\
\\
{\bf Example 4.1} ([D1; $\S$ 2], [D4; $\S$ 2]). Let $\beta \in \mathbb{C}^{\times}$ and  $z \in \mathbb{C}$. \\
$$ \underline{\prod}_{\nu = 0}^{\infty}  \beta(z + \nu) = \beta^{1/2 - z} \left(  \frac{1}{2\pi} \Gamma(z) \right)^{-1},$$
where $\Gamma(z)$ is the gamma function.
$$ \underline{\prod}_{\nu \in \mathbb{Z}} \beta(z + \nu) = \left\{  \begin{array}{ll} 1 - e^{-2\pi\sqrt{-1} z} &  \; \mbox{if}\; {\rm Im}(\beta) >0,\\
1 - e^{2\pi\sqrt{-1} z} &  \; \mbox{if}\; {\rm Im}(\beta) <0. \\
\end{array} \right.$$
\\
Using the formulas in Example 4.1 and the classical Lefschetz trace formula, we can show the regularized determinant formula for the zeta function of a Riemannian foliated dynamical system of type (i) as follows.\\
\\
{\bf Theorem 4.2.} {\it Let $\frak{S} = (M, g, {\cal F}, \phi)$ be a 3-dimensional Riemannian foliated dynamical system of type {\rm (i)} so that $M$ is a fiber bundle over $S^1$ whose typical fiber is an oriented, connected, closed Riemannian surface $S$, ${\cal F}$ is the bundle foliation, and $\phi$ is leafwise homotopic to the suspension flow of the monodromy diffeomorphism $f : S \rightarrow S$. As in Example 3.3, we assume that $f$ is an Anosov diffeomorphism when $S$ is a torus and that $f$ is an Axiom A diffeomorphism when the genus of $S$ $\geq 2$. Let $\Theta$ denote the infinitesimal operator of $\phi^t$.}\\
(1) {\it  We have}
$$ {\rm Sp}(\Theta \, | \, \bar{H}_{\cal F}^n(M)_{\mathbb{C}}) = \{ \log \alpha + 2\pi\sqrt{-1} \nu \; | \; \alpha \in {\rm Sp}(H^n(f)\,|\,H^n(S, \mathbb{C})), \nu \in \mathbb{Z}  \}.$$
(2) {\it The regularized determinants ${\rm det}_{\infty}(\Theta \, | \, \bar{H}_{\cal F}^n(M)_{\mathbb{C}})$ are defined and we have the following formula}: 
$$ \zeta_{\frak{S}}(s) = \prod_{n=0}^2 {\rm det}_{\infty}( {\rm id} - \Theta \; | \; \bar{H}_{\cal F}^n(M)_{\mathbb{C}})^{(-1)^{n+1}}.$$
\\
{\it Proof.} (1) By Example 1.6, $\bar{H}^n_{\cal F}(M)_{\mathbb{C}}$ is isomorphic to the smooth sections of a flat vector bundle ${\cal E}$ over $S^1$ whose typical fiber is $H^n(S, \mathbb{C})$, defined by the suspension of $H^n(f)$. Let $\alpha \in  {\rm Sp}(H^n(f)\; |\; H^n(S,\mathbb{C}))$ so that there is $u \in H^n(S, \mathbb{C})$ such that $H^n(f)(u) = \alpha u$. Then $u$ defines a constant section of ${\cal E}$ and so an element $\bar{u} \in \bar{H}^n_{\cal F}(M)$ satisfying 
$$ \phi^{t *} (\bar{u}) = \alpha^t \cdot \bar{u},$$
since we take $S^1 = \mathbb{R}/\mathbb{Z}$, and therefore 
$$ \Theta(\bar{u}) = \log \alpha \cdot \bar{u}.$$ 
Now let $v_{\nu} = e^{2\pi \sqrt{-1} \nu s} \cdot \bar{u} \in \bar{H}^{n}_{\cal F}(M)_{\mathbb{C}}$ for $s \in S^1 = \mathbb{R}/\mathbb{Z}$, $\nu \in \mathbb{Z}$. The functions $v_{\nu}$  $(\nu \in \mathbb{Z})$ form a complete orthonormal system of eigenfuctions of the Laplacian $\Delta_{S^1} = - d^2/ds^2$, with corresponding eigenvalues
$4\pi \nu^2$. So, if $H^n(S,\mathbb{C})$ has orthonormal frame of eigenvectors $u_j$ of $H^n(f)$ with eigenvalues $\alpha_j$, then the elements $e^{2\pi \sqrt{-1} \nu s} \cdot \bar{u}_j$ of $\bar{H}_{\cal F}^n(M)_{\mathbb{C}}$  form a complete orthonormal system of $\bar{H}_{\cal F}^n(M)_{\mathbb{C}}$. Moreover we have
$$ \Theta(e^{2\pi \sqrt{-1}\nu s} \cdot \bar{u}_j) = (2\pi \sqrt{-1} \nu + \log \alpha) e^{2\pi \sqrt{-1} \nu} \bar{u}_j.$$
If there is not a complete orthonormal system of $H^n(S,\mathbb{C})$ consisting of eigenvectors of $H^n(f)$, then we get all eigenvalues of $\Theta$ by using all eigenvalues $\alpha_j$ of $H^n(f)$ on $H^n(S,\mathbb{C})$ with the above procedure. \\
(2) By (1), we have
$$ \begin{array}{l} \displaystyle{\prod_{n=0}^2 {\rm det}_{\infty}(s\; {\rm id} - \Theta \; | \; \bar{H}_{\cal F}^n(M)_{\mathbb{C}})^{(-1)^{n+1}}} \\
       = \displaystyle{ \prod_{n=0}^2  \prod_{\alpha \in {\rm Sp}(H^n(f))} \underline{\prod}_{\nu \in \mathbb{Z}} ( s - (\log \alpha + 2\pi \sqrt{-1} \nu))^{(-1)^{n+1}} }\\
              = \displaystyle{  \prod_{n=0}^2  \prod_{\alpha \in {\rm Sp}(H^n(f))} \underline{\prod}_{\nu \in \mathbb{Z}} 2\pi\sqrt{-1} \left( \frac{s - \log \alpha}{2\pi\sqrt{-1}} + \nu \right)^{(-1)^{n+1}} }\\
     = \displaystyle{ \prod_{n=0}^2  \prod_{\alpha \in {\rm Sp}(H^n(f))} (1- \alpha e^{-s})^{(-1)^{n+1}}} \;\;\;\; (\mbox{Example 4.1})\\
     = \displaystyle{ \prod_{n=0}^2 {\rm det}(  {\rm id} - e^{-s} H^n(f) \; | \; H^n(S,\mathbb{C}))^{(-1)^{n+1}} }\\
     = \displaystyle{ \exp \left( \sum_{n=0}^2 (-1)^n  \sum_{k=1}^{\infty} \frac{1}{k} e^{-sk} {\rm tr}(H^n(f^k))  \right). } 
\end{array} \leqno{(4.2.1)}$$
The classical Lefschetz trace formula tells us the equality
$$ \sum_{p \in {\rm Fix}(f^k)} \varepsilon_p = \sum_{n=0}^2 (-1)^n  {\rm tr}(H^1(f^k)), \leqno{(4.2.2)} $$ 
where ${\rm Fix}(f^k) := \{ p \in S \, | \, f^k(p) = p\}$ and  $\varepsilon_p := {\rm sgn} \det ({\rm id} - T_p(f^k))$.  In fact, for every $p \in \operatorname{Fix}(f^k)$,  either $\epsilon_p = -1$ for all $k$, or $\epsilon_p = -1$ for all $k$, as we have seen in Example 3.3. 
By the standard computation, we have the equality 
$$ \exp \left(  \sum_{k=1}^{\infty} \frac{1}{k} \left(\sum_{p \in {\rm Fix}(f^k)} \varepsilon_p  \right) e^{-ks} \right) = \prod_{\frak{p} \in {\cal P}^f} (1 - e^{-s \#(\frak{p})})^{-\varepsilon_{\frak{p}}}, \leqno{(4.2.3)}$$
 where ${\cal P}^f$ denotes the set of finite orbits of the iteration of $f$ on $S$ and $\varepsilon_{\frak{p}} = \varepsilon_p$ for the orbit $\frak{p}$ containing $p$. By (4.2.1), (4.2.2) and (4.2.3), we obtain
$$ \prod_{n=0}^2 {\rm det}_{\infty}( s\; {\rm id} - \Theta \; | \; \bar{H}_{\cal F}^n(M))^{(-1)^{n+1}} = \prod_{\frak{p} \in {\cal P}^f} (1 - e^{-s \#(\frak{p})})^{-\varepsilon_{\frak{p}}}.
\leqno{(4.2.4)} $$
Since the correspondence $\gamma \mapsto \gamma \cap S$ gives the bijection ${\cal P} \stackrel{\sim}{\rightarrow} {\cal P}^f$ such that $\varepsilon_{p} = \varepsilon_{\gamma}$ for $p \in S \cap \gamma$ and $\#(\frak{p}) = \ell(\gamma)$, we have
$$ \prod_{\frak{p} \in {\cal P}^f} (1 - e^{-s \#(\frak{p})})^{-\varepsilon_{\frak{p}}} = \zeta_{\frak{S}}(s).  \leqno{(4.2.5)}$$
By (4.2.4) and (4.2.5), we obtain
$$ \prod_{n=0}^2 {\rm det}_{\infty}(s\; {\rm id} - \Theta \; | \; \bar{H}_{\cal F}^n(M))^{(-1)^{n+1}} = \zeta_{\frak{S}}(s). \;\; \Box $$
\\

\begin{center}
{\bf 5. Examples of 3-dimensional Riemannian foliated dynamical systems of type (ii)}
\end{center}

In this section, we provide examples of 3-dimensional Riemannian foliated dynamical systems $\frak{S} = (M, g, {\cal F}, \phi)$ of type (ii) $(\mbox{cf. Theorem 1.3})$, which satisfy the following  conditions:\\
(A1) $\phi$ is simple. \\
(A2) The condition (4) in Lemma 3.1, namely, $\displaystyle{ a := \limsup_{x \rightarrow \infty} \frac{ \log \nu(x)}{x} }$ exists.\\
(A3) $\varepsilon_{\gamma}(k) = \varepsilon_{\gamma}$ for any closed orbit $\gamma$ and any $k \in \mathbb{N}$.\\
(A4) $\phi$ is leafwise homotopic to a foliated flow $\psi$ on $(M, {\cal F})$, which is isometric.\\
\\
Note that the conditions (A1)$\sim$(A3) are satisfied for a 3-dimensional Riemannian foliated dynamical system of type (i), as was shown in Example 3.3 and Theorem 4.2.
But (A4) may not be satisfied by systems of type (i).\\

As for (A4), we have the following equivalent condition, which is convenient for use.\\
\\
{\bf Lemma 5.1.} {\it Let $\phi$ and $\psi$ be foliated flows on $(M, {\cal F})$. Let $X$ and $Y$ be the vector fields $($the infinitesimal generators$)$ attached to $\phi$ and $\psi$, respectively.
Then $\phi$ is leafwise homotopic to $\psi$  if and only if $Y - X$ is tangent to leaves, namely, the sections of the normal bundle $T(M)/T({\cal F})$ defined by $X$ and $Y$ are same.}\\
\\
{\it Proof.} The smooth family of foliated flows $H_s$ over $s \in [0,1]$ with infinitesimal generator $(1-s)X + sY$ give the leafwise homotopy connecting $\phi$ and $\psi$. $\;\; \Box$\\
\\
{\bf Example 5.2.} Let $S$ be an oriented, connected, closed surface of genus at least 2. Take a homomorphism $\rho : \pi_1(S) \to \mathbb{R}$ such that ${\rm Im}(\rho)$ is dense in $\mathbb{R}$. Let $\tilde{S} \to S$ be the universal covering so that $\pi_1(S)$  acts canonically on $\tilde{S}$ from the right. Then one has the right action of $\pi_1(S)$ on $\tilde{S} \times \mathbb{R}$ defined by 
$$(\tilde{p}, r).\alpha := (\tilde{p}.\alpha, r - \rho(\alpha))$$
 for $\tilde{p} \in \tilde{S}, r \in \mathbb{R}$ and $\alpha \in \pi_1(S)$, and we define  $M$ by the quotient
$$ M := (\tilde{S} \times \mathbb{R})/\pi_1(S).$$
The foliation on $\tilde{S} \times \mathbb{R}$ with leaves $\tilde{S} \times \{ r \}$ $(r \in \mathbb{R})$ induces the foliation ${\cal F} = \{ L_r \}_{r \in \mathbb{R}}$ on $M$, where $L_r$ is the image of $\tilde{S} \times \{ r \}$ in $M$. Since $\rho$ has the dense image, the leaves $L_r$ are dense in $M$. Let $\pi : M \rightarrow S$ be the projection induced by the 1st projection $\tilde{S} \times \mathbb{R} \to \tilde{S}$ and the covering $\tilde{S} \to S$. Then $\pi$ is a fiber bundle with typical fiber $S^1$ and the horizontal leaves $L_r$ are transverse to the fibers of $\pi$. 
Take a hyperbolic metric on $S$ and take a bundle-like metric $g$ on $M$ so that $\pi$ is a Riemannian submersion.  We can assume the fibers of $\pi$ are of length one. 
Let $\psi$ be the vertical flow on $M$ defined by
$$ \psi^t([\tilde{p}, r]) := [\tilde{p}, t + r],$$
which acts isometrically on the leafwise differential forms on $M$. The orbits of $\psi$ are the fibers of $\pi$, and therefore they are closed of length one. Let $X$ be the vector field attached to $\psi^t$. Take a vector field $\bar{Y}$ on $S$ and let $Y$ be the lift of $\bar{Y}$ to $M$, which is tangent to the leaves $L_r$ of ${\cal F}$ (the horizontal lift of $\bar{Y}$). Let $Z := X + Y$ and let $\phi$ be the flow of $Z$, which is an ${\cal F}$-foliated flow. By Lemma 5.1, $\phi$ and $\psi$ are leafwise homotopic. 
Thus we obtain a $3$-dimensional Riemannian foliated dynamical system $\frak{S} = (M,g,{\cal F},\phi)$ of type (ii) satisfying (A4). 

As $\bar{Y}$, consider the gradient vector field of a Morse function $h$ on $S$, $\bar{Y} = {\rm grad}(h)$. Then any closed orbit of $\phi$ is simple. Moreover, the length of closed orbits is one, and the closed orbits  correspond to the zero points of $\bar{Y}$. If $\bar{Y}$ has $N$ zero points, the number of closed orbits of $\phi$ is $N$. Therefore the conditions (A1), (A2) are satisfied for $\frak{S} = (M,g,{\cal F},\phi)$. 

Lastly, we shall see that (A3) is satisfied. Let $\bar\phi^t$ be the flow of $\bar{Y}$. Note that $X$ is vertical and so the projection $\pi: M \to S$ is equivariant with respect to the flows $\phi^t$ on $M$ and $\bar{\phi}^t$ on $S$. Let $\gamma$ be a closed orbit and let $x \in \gamma$. By the equivariance, $\gamma$ corresponds to a zero point $\pi(x)$ of $\bar{Y}$, which is a fixed point of $\bar{\phi}^t$. Moreover, for every $k\in\mathbb N$, the tangent map $T_x(\phi^k) : T_x({\cal F}) \to T_x({\cal F})$ corresponds to the tangent map $T_{\pi(x)}(\bar{\phi}^k) : T_{\pi(x)}(S) \to T_{\pi(x)}(S)$. So each index $\varepsilon_{\gamma}(k)$ is equal to the sign of the determinant of $T_{\pi(x)}(\bar{\phi}^{k}) - {\rm id} : T_{\pi(x)}(S) \to T_{\pi(x)}(S)$. 
Taking Morse coordinates, we have 3 local models:\\
\\
$\cdot$ Case of Morse index 0: \\
$h(x,y) = -x^2-y^2.$ \\
$\bar{Y} = \nabla h = (-2x, -2y),$ $\; \bar\phi^t(x,y) = (e^{-2t}, e^{-2t}).$\\
$T_{(0,0)}(\bar\phi^t) \equiv \begin{pmatrix} -2 & 0 \\ 0 & -2 \end{pmatrix},$ $T_{(0,0)}(\bar\phi^t) -{\rm id} \equiv \begin{pmatrix} -3 & 0 \\ 0 & -3 \end{pmatrix}$ for any $t$.\\
Therefore $\varepsilon_\gamma(k) = \varepsilon_{\gamma} = 1$ for any $k \in \mathbb{N}$ and $\gamma$.
\\
$\cdot$ Case of Morse index 1: \\
$h(x,y) = x^2-y^2.$\\
$\bar Y = \nabla h = (2x, -2y),$ $\; \bar\phi^t(x,y) = (e^{2t}, e^{-2t}).$\\
$T_{(0,0)}(\bar\phi^t) \equiv \begin{pmatrix} 2 & 0 \\ 0 & -2 \end{pmatrix},$  $T_{(0,0)}(\bar\phi^t) -{\rm id} \equiv \begin{pmatrix} 1 & 0 \\ 0 & -3 \end{pmatrix}$ for any $t$. \\
Therefore $\varepsilon_\gamma(k) = \varepsilon_{\gamma} = -1$ for any $k \in \mathbb{N}$ and $\gamma$.\\
$\cdot$ Case of Morse index 2: \\
$h(x,y) = x^2+y^2. $\\
$\bar Y = \nabla h = (2x, 2y),$ $\; \bar\psi^t(x,y) = (e^{2t}, e^{2t}).$\\
$T_{(0,0)}(\bar\phi^t) \equiv \begin{pmatrix} 2 & 0 \\ 0 & 2 \end{pmatrix},$ $T_{(0,0)}(\bar\phi^t) -{\rm id} \equiv \begin{pmatrix} 1 & 0 \\ 0 & 1 \end{pmatrix}$ for any $t$.\\
Therefore $\varepsilon_\gamma(k) = \varepsilon_{\gamma} = 1$ for any $k \in \mathbb{N}$ and $\gamma$.\\
\\
So $\varepsilon_\gamma(k) = \varepsilon_{\gamma}$ depends on the Morse index in the above way. In particular, the condition (A3) is satisfied. We can group the closed orbits according to the Morse index, obtaining that the zeta function $\zeta_{\frak{S}}(s)$ has 3 factors, 2 of them with exponent -1, and one of them with exponent 1.

Finally, we show that the reduced leafwise cohomology $\bar{H}_{\cal F}^1(M)$ is of infinite dimension for this $\frak{S}$. Since the genus of $S$ is al least 2, there are two transverse circles $C$ and $C'$ in $S$ such that $[C], [C'] \in {\rm Ker}(\rho), [C], [C'] \neq 1$ and $C$ and $C'$ have non-trivial intersection number and such that the restriction of ${\cal F}$ to $\pi^{-1}(C)$ and $\pi^{-1}(C')$ are product foliations of dimension one. This can be achieved when $S$ is a surface of genus $\ge 2$. Then the foliation of $\pi^{-1}(C)$ has infinitely many transverse invariant measures $\mu_k$ (invariant distributions are also valid) with disjoint support, and the foliation of $\pi^{-1}(C')$ has infinitely many smooth basic functions $f_l$ with disjoint support. We can also assume that the restrictions of $\mu_k$ and $f_l$ to any common transversal $\pi^{-1}(p) \; (p \in C \cap C')$ satisfy $\mu_k(f_l) = \delta_{kl}.$ Take also a cohomology class in $S$ which is Poincar\'{e}  dual of $C$ and represented by a closed 1-form $\eta$ supported in a small tubular neighborhood $V$ of $C'$. It is also true that the restriction of ${\cal F}$ to $\pi^{-1}(V)$ is a product foliation. Then the 1-forms $f_l\pi^*\eta$ are leafwise closed, defining leafwise reduced cohomology classes of ${\cal F}$, and their restrictions to $\pi^{-1}(C)$ satisfy 
$\langle [\mu_k], [f_l\pi^*\eta] \rangle \ne 0 \; \mbox{if and only if} \; k=l.$
Here, we consider the transverse invariant measures $\mu_k$ as foliation cycles of the foliation on $\pi^{-1}(C)$ that can be paired with its leafwise reduced cohomology classes of top degree on $\pi^{-1}(C)$, and the leafwise reduced cohomology classes $[f_l\pi^*\eta]$ are restricted to $\pi^{-1}(C)$ in this pairing. Hence the leafwise reduced cohomology $\bar{H}_{\cal F}^1(M)$ is of infinite dimension.\\

\begin{center}
{\bf 6. Regularized determinant formulas for the zeta functions for certain 3-dimensional Riemannian foliated dynamical systems of type (ii)}
\end{center}

In this section, we prove a regularized determinant formula of the dynamical zeta function for a Riemannian foliated dynamical system $\frak{S} = (M, g, {\cal F}, \phi)$ of type (ii), which  satisfies the conditions (A1) $\sim$ (A4) in the section 5. We write ${\rm Sp}^n(\Theta)$ for ${\rm Sp}(\Theta | \bar{H}^n_{\cal F}(M)_{\mathbb{C}})$. \\

For $n = 0, 1, 2$, we define the {\it dynamical spectral $\xi$-function} $\xi_n(z,s)$ of $\frak{S}$  by
$$ \xi_n(z,s) := \sum_{\alpha \in {\rm Sp}^n(\Theta)} (s - \alpha)^{-z},$$
where ${\rm Re}(z), {\rm Re}(s) > 1$, ${\rm Arg}(s - \rho) \in (-\pi/2, \pi/2)$. By Theorem 1.8, we have $\xi_0(z, s) = \xi_2(z,s) = s^{-z}$. 
In the following, we relate the alternative sum 
$$\displaystyle{ \sum_{n=0}^2 (-1)^n \xi_n(z,s)  = 2 s^{-z} - \sum_{\rho \in {\rm Sp}^1(\Theta)} (s - \rho)^{-z}} \leqno{(6.1)}$$
with the zeta function $\zeta_{\frak{S}}(s)$, by applying Theorem 2.5 (3).

Note  by (A1) and Lemma 2.2 that the set ${\cal P}$ of closed orbits is countable. Moreover $\mu$ is a Radon measure on $\mathbb R_+$ by (1) of Proposition 3.4 and (A2). For $z,s \in \mathbb C$, consider the function $F : \mathbb{R}_+ \rightarrow \mathbb{C}$ defined by
$$F(x) :=  x^{z-1}e^{-sx}.$$
For every $w\in\mathbb C$, let $\mu_w$ be the measure on $\mathbb R_+$ defined by the function $e^{wx}$. By (A2) and (2) of Proposition 3.4,  if $\operatorname{Re}(s) > a$ ($a$ is given in (A2)), then the function $F$ is absolutely integrable with respect to $|\mu|$. If moreover $\operatorname{Re}(w) \le a$ and $\operatorname{Re}(z) \ge 1$, then $F$ is also absolutely integrable with respect to $|\mu_w|$.
\\
\\
{\bf Proposition 6.2.} {\it If $\operatorname{Re}(z)>5$, then we have
\[
2\mu_0(F) + \sum_{\rho\in\operatorname{Spec}^1(\Theta)} \mu_\rho(F) = \mu(F),
\]
where the series in the left-hand side of the equality is absolutely convergent.}
\\
\\
{\it Proof.} First, let us prove that the series in the left-hand side of the stated equality is absolutely convergent. We have
\begin{align*}
F'(x) &= \big((z-1)x^{z-2} - sx^{z-1}\big) e^{-sx}, \\
F''(x) &= \big((z-1)(z-2)x^{z-3} - 2s(z-1)x^{z-2} + s^2x^{z-1}\big) e^{-sx}, \\
F'''(x) &= \big((z-1)(z-2)(z-3)x^{z-4} - 3s(z-1)(z-2)x^{z-3} \\
&\phantom{={}} {} + 3s^2(z-1)x^{z-2} - s^3x^{z-1}\big) e^{-sx}, \\
F^{(4)}(x) &= \big((z-1)(z-2)(z-3)(z-4)x^{z-5} - 4s(z-1)(z-2)(z-3)x^{z-4} \\
&\phantom{={}} {} + 6s^2(z-1)(z-2)x^{z-3} - 4s^3(z-1)x^{z-2} + s^4x^{z-1}\big) e^{-sx}.
\end{align*}
  Since ${\rm Re}(z) > 5$, we see that, for $0\le k \le 4$,  $F^{(k)}$ is defined and vanishes at $x=0$ and  $F^{(k)}(x) \to 0$ as $x\to\infty$. For all $\rho\in \operatorname{Spec}^1(\Theta)$,  we have $\operatorname{Re}(\rho) = 0$  by Theorem 2.5 (1), and therefore $F$ is absolutely integrable with respect to $|\mu_\rho|$. So, using integration by parts four times, we obtain
\[
\mu_\rho(F)=\int_0^\infty F(x)e^{\rho x}\,dx = \rho^{-4} \int_0^\infty F^{(4)}(x)\,e^{\rho x}\,dx = \rho^{-4}\,\mu_\rho(F^{(k)})\;.
\]
Moreover, there is some $c_1>0$, depending only on $s$, such that $|\mu_{\rho}(F^{(k)})| < c_1$ for all $\rho\in \operatorname{Spec}^1(\Theta)$, obtaining
\[
|\mu_\rho(F)| \le c_1 |\rho|^{-4}.
\]
Then, by (2) of Theorem 2.5, 
\[
\sum_{\rho\in\operatorname{Spec}^1(\Theta)} |\mu_\rho(F)| = \sum_{j=0}^\infty |\rho_j|^{-4}|\mu_{\rho_j}(F^{(4)})| \le c_2 \sum_{j=0}^\infty j^{-4/3} <\infty,
\]
for some constant $c_2 >0$. 

Next, let us prove that the stated equality holds. Consider a sequence $f_i$ in $C^\infty_{\text{\rm c}}(\mathbb R)$ ($i=1,2,\dots$) such that $|f_i| \le |F|$ and $f_i = F$ on $[1/i, i]$. By the properties mentioned before the statement  and Lebesgue's dominated convergence theorem, we have $\mu(f_i) \to \mu(F)$ and $\mu_\rho(f_i) \to \mu_\rho(F)$ as $i\to\infty$ for all $\rho\in\operatorname{Spec}^1(\Theta)$. Since the equality of (3) in Theorem 2.5 holds in the space of distributions, we get
\[
2\mu_0(f_i) + \sum_{\rho\in\operatorname{Spec}^1(\Theta)} \mu_\rho(f_i) = \mu(f_i).
\]
Using the absolute convergence of the series on the left-hand side of the stated equality and the absolute integrability  of $F$ with respect to every $|\mu_\rho|$ on $[0,\infty)$, it follows from Lebesgue's dominated convergence theorem that
$$ \begin{array}{ll} 
\displaystyle{ 2\mu_0(F) + \sum_{\rho\in\operatorname{Spec}^1(\Theta)} \mu_\rho(F) } & = \displaystyle{ \lim_{i\to\infty}\Big( 2\mu_0(f_i) + \sum_{\rho\in\operatorname{Spec}^1(\Theta)} \mu_\rho(f_i) \Big) }\\
&= \displaystyle{ \lim_{i\to\infty} \mu(f_i) } \\
& = \mu(F). \;\; \Box
\end{array}
$$
\\
Assume $\operatorname{Re}(z)>5$.
Proposition 6.2 yields the equality
$$ 2\mu_0(F) -  \sum_{\rho \in {\rm Sp}^1(\Theta)} \mu_{\rho}(F)  = \sum_{\gamma \in {\cal P}} \ell(\gamma) \sum_{k = 1}^{\infty} \varepsilon_{\gamma}(k) F(k\ell(\gamma)). \leqno{(6.2)} $$
By [E; 1.5.1 (34)], we have 
$$ \mu_w(F) = \Gamma(z)(s-w)^{-z} \;\;  \mbox{for}\; {\rm Re}(s) > {\rm Re}(w), \leqno{(6.3)}$$
where ${\rm Arg}(s-w) \in (-\pi/2, \pi/2)$.  Therefore, by (6.1), (6.2) and (6.3), we have 
$$ \sum_{n=0}^2 (-1)^n \xi_n(z,s) = \frac{1}{\Gamma(z)} \sum_{\gamma \in {\cal P}} \ell(\gamma) \sum_{k = 1}^{\infty} \varepsilon_{\gamma}(k) F(k\ell(\gamma)) \leqno{(6.4)}$$
for ${\rm Re}(z) > 5$ and ${\rm Re}(s) > {\rm max}\{ 1, a \}$. Let $\mathcal{C}$ be the contour consisting of the lower edge of the cut from $-\infty$ to $-\delta$, the circle $x = \delta e^{\sqrt{-1} \varphi}$ with $-\pi \leq \varphi \leq \pi$ and the upper edge of the cut from $-\delta$ to $- \infty$. By [E; 1.6 (6)], we have, for $c >0$,
$$ \frac{c^{z-1}}{\Gamma(z)} = \frac{1}{2\pi \sqrt{-1}} \int_{\cal C} t^{-z} e^{ct} dt. $$
Therefore, for $0 < \delta < {\rm Re}(s) - 1$, we have
$$ \begin{array}{l} \displaystyle{ \frac{1}{\Gamma(z)}  \sum_{\gamma \in {\cal P}} \ell(\gamma) \sum_{k = 1}^{\infty} \varepsilon_{\gamma}(k) F(k\ell(\gamma))} \\ 
=   \displaystyle{ \frac{1}{\Gamma(z)} \sum_{\gamma \in {\cal P}} \ell(\gamma) \sum_{k = 1}^{\infty} \varepsilon_{\gamma}(k) (k \ell(\gamma))^{z-1} e^{-s k \ell(\gamma)} }\\
 =  \displaystyle{ \sum_{\gamma \in {\cal P}} \ell(\gamma)  \sum_{k = 1}^{\infty} \varepsilon_{\gamma}(k) e^{-s k \ell(\gamma)} \frac{  (k \ell(\gamma))^{z-1}  }{\Gamma(z)} }\\
 =  \displaystyle{ \frac{1}{2\pi \sqrt{-1}} \sum_{\gamma \in {\cal P}} \ell(\gamma)  \sum_{k = 1}^{\infty} \varepsilon_{\gamma}(k) \int_{\cal C} t^{-z} e^{k \ell(\gamma) (t-s)} dt }
\end{array} \leqno{(6.5)}
$$
By (6.4) and (6.5), we obtain the following\\
\\
{\bf Proposition 6.6.} {\it Notations being as above, we have }
$$ \sum_{n=0}^2 (-1)^n \xi_n(z,s)  =   \displaystyle{ \frac{1}{2\pi \sqrt{-1}} \sum_{\gamma \in {\cal P}} \ell(\gamma)  \sum_{k = 1}^{\infty} \varepsilon_{\gamma}(k) \int_{\cal C} t^{-z} e^{k \ell(\gamma) (t-s)} dt }$$
{\it for $\operatorname{Re}(z)>5$ and ${\rm Re}(s) > {\rm max}\{1, a \}$.}\\
\\
By (A2), Lemma 3.1 and Proposition 3.2,  $\zeta_{\frak{S}}(s)$ converges absolutely for ${\rm Re}(s) > {\rm max}\{ a, \frac{\log 2}{m} \}$, where $m := {\rm Inf}\{ \ell(\gamma) \; | \; \gamma \in {\cal P} \}$.
Then we have
$$ \frac{\zeta_{\frak{S}}'}{\zeta_{\frak{S}}}(s) = \partial_s \log \zeta_{\frak{S}}(s) = - \sum_{\gamma \in {\cal P}} \ell(\gamma) \varepsilon_{\gamma} \sum_{k=1}^{\infty} e^{-s k \ell(\gamma)} \;\;\;\; ({\rm Re}(s) > r).\leqno{(6.7)} $$
We set
$$ b := {\rm max}\{ 1, a, \frac{\log 2}{m} \}.$$
By  Proposition 6.6, (6.7)  and (A3), we obtain the following theorem, which relates the dynamical $\xi$-functions with the zeta function.
\\
\\
{\bf Theorem 6.9.} {\it For ${\rm Re}(s) > b$ and $\operatorname{Re}(z)>5$, we have}
$$  \sum_{n=0}^2 (-1)^n \xi_n(z,s) =  \displaystyle{- \frac{1}{2\pi \sqrt{-1}} \int_{\cal C} t^{-z} \frac{\zeta_{\frak{S}}'}{\zeta_{\frak{S}}}(s-t) dt.} 
$$ 
To get Theorem 6.9,  we interchanged the infinite sum and the integral by Lebesgue dominated converge: observe that $\displaystyle{ \int_{-\infty}^{-\delta} t^{-{\rm Re}(z)} \frac{\zeta_{\frak{S}}'}{\zeta_{\frak{S}}}({\rm Re}(s) - t) dt}$ converges absolutely. The integral in Theorem 6.9 defines a holomorphic function of $z$. So Theorem 6.9 yields \\
\\
{\bf Corollary 6.10.} For ${\rm Re}(s) > b$ and $n = 0, 1,2$, there is a holomorphic continuation of $z \mapsto \xi_n(z,s)$ to $ \mathbb{C} \setminus \{ 1 \}$, also denoted by $\xi_n(z,s)$.\\

By the definition of the regularized determinant, we have 
$$ {\rm det}_{\infty} ( s\; {\rm id} - \Theta\; | \; \bar{H}^n_{\cal F}(M)_{\mathbb{C}}) = \exp( - \partial_z \xi_n(z,s)|_{z = 0}). $$
By Theorem 6.9, we have, for ${\rm Re}(s) > b$, 
$$ \prod_{n=0}^2 {\rm det}_{\infty} ( s \;{\rm id} - \Theta\; | \; \bar{H}^n_{\cal F}(M)_{\mathbb{C}})^{(-1)^{n+1}}
 = \exp((\partial_z \displaystyle{\frac{1}{2\pi \sqrt{-1}} \int_{\cal C} t^{-z} \frac{\zeta_{\frak{S}}'}{\zeta_{\frak{S}}}(s-t) dt})|_{z=0} ). \leqno{(6.11)}$$
Here we have
$$ \begin{array}{ll}
 ( \partial_z \displaystyle{\frac{1}{2\pi \sqrt{-1}} \int_{\cal C} t^{-z} \frac{\zeta_{\frak{S}}'}{\zeta_{\frak{S}}}(s-t) dt})|_{z=0} & 
=  \displaystyle{\frac{1}{2\pi \sqrt{-1}} \int_{\cal C} (\partial_z t^{-z} \frac{\zeta_{\frak{S}}'}{\zeta_{\frak{S}}}(s-t))|_{z=0} dt} \\
& =  \displaystyle{\frac{1}{2\pi \sqrt{-1}} \int_{\cal C} \frac{\zeta_{\frak{S}}'}{\zeta_{\frak{S}}}(s-t) \log (t) dt} \\
&  =  \displaystyle{\frac{1}{2\pi \sqrt{-1}} \int_{-\infty}^0 \frac{\zeta_{\frak{S}}'}{\zeta_{\frak{S}}}(s-t) (\log |t| - \sqrt{-1}\pi) dt} \\
& \;\;\;\; + \displaystyle{\frac{1}{2\pi \sqrt{-1}} \int_{0}^{-\infty} \frac{\zeta_{\frak{S}}'}{\zeta_{\frak{S}}}(s-t) (\log |t| + \sqrt{-1}\pi)  dt}\\
& = \displaystyle{ \int_0^{-\infty} \frac{\zeta_{\frak{S}}'}{\zeta_{\frak{S}}}(s-t) dt} \\
& = \displaystyle{ - \int_0^{+\infty} \frac{\zeta_{\frak{S}}'}{\zeta_{\frak{S}}}(s+t) dt} \\
& =  \displaystyle{ \sum_{\gamma \in {\cal P}} \sum_{k=1}^{\infty} \varepsilon_{\gamma} \frac{e^{-s k \ell(\gamma)}}{k}  }\\
& = \log \zeta_{\frak{S}}(s). 
\end{array} \leqno{(6.12)}
$$
Hence, by (6.11) and (6.12), we obtain\\
\\
{\bf Theorem 6.13.} {\it Notations being as above, we have
$$ \zeta_{\frak{S}}(s) = \prod_{n=0}^2 {\rm det}_{\infty} ( s\; {\rm id} - \Theta\; | \; \bar{H}^i_{\cal F}(M)_{\mathbb{C}})^{(-1)^{n+1}} $$
for ${\rm Re}(s) > b$.}
\\

\begin{flushleft}
{\bf References}\\
{[\'{A}H]} J.A. \'{A}lvarez L\'{o}pez, G. Hector, The dimension of the leafwise reduced cohomology, Amer. J. Math., {\bf 123}, no. 4, 2001, 607--64.\\
{[\'{A}K1]} J. A. \'{A}lvarez L\'{o}pez, Y. A. Kordyukov, Long time behaviour of
leafwise heat flow for Riemannian foliations, Compositio Mathematica, {\bf 125},
no. 2, 2001, 129--153.\\
{[\'{A}K2]} J. A. \'{A}lvarez L\'{o}pez, Y. A. Kordyukov, Distributional Betti
numbers of transitive foliations of codimension one, In: Foliations: Geometry
and Dynamics, World Scientific, 2002, 159--183.\\
{[\'{A}KL1]} J. A. \'{A}lvarez L\'{o}pez, Y. A. Kordyukov, E. Leichtman,  Simple foliated flows, Tohoku Mathematical Journal, {\bf 74}, no. 1, 2022, 53--81.\\
{[\'{A}KL2]} J. A. \'{A}lvarez L\'{o}pez, Y. A. Kordyukov, E. Leichtman,  A trace formula for foliated flows, arXiv:2402.06671\\
{[BGV]} N. Berline, E. Getzler. and M. Vergne, Heat kernels and Dirac operators. Springer Verlag, New York, 1992.\\
{[Bo]} R. Bowen, Topological entropy and axiom A, Proc. Symp. in Pure Math. AMS, {\bf 14}, Global Analysis (1970), 23--41.\\
{[C]} A. Connes, Trace formula in noncommutative geometry and the zeros of the Riemann zeta function, Selecta Mathematica, New series {\bf 5}, 1999, 29--106.\\
{[CC]} A. Connes, K. Consani, Cyclic homology, Serre's local factors and $\lambda$-operations, J. K-Theory, {\bf 14}, no. 1, 2014, 1--45.\\
{[D1]} C. Deninger, On the $\Gamma$-factors attached to motives, Inventiones Mathematicae, {\bf 104}, 1991, 245--261.\\
{[D2]} C. Deninger, Local $L$-factors of motives and regularized determinants, Inventiones Mathematicae, {\bf 107}, 1992, 135--150.\\
{[D3]} C. Deninger, Lefschetz trace formulas and explicit formulas in analytic number theory, Journal f\"{u}r Reine und Angewandte Mathematik, {\bf 441}, (1993), 1--15.\\
{[D4]} C. Deninger, Motivic $L$-functions and regularized determinants, Proc. Symp. Pure Math., AMS, {\bf 55}, 1, 1994, 707--743.\\
{[D5]} C. Deninger, Evidence for a cohomological approach to analytic number theory,  In: Joseph, A., Mignot, F., Murat, F., Prum, B., Rentschler, R. (eds), First European Congress of Mathematics, Progress in  Mathematics {\bf 3}, Birkh\"{a}user Boston, Boston, MA, 491--510.\\
{[D6]} C. Deninger, Some analogies between number theory and dynamical 
systems on foliated spaces, Documenta Mathematica, Jahresbericht der Deutschen Mathematiker-Vereinigung, Extra Volume International Congress of Mathematicians I,  
1998, 23--46.\\
{[D7]} C. Deninger, On dynamical systems and their possible significance for
arithmetic geometry, In: Reznikov, A., Schappacher, N. (eds), Regulators in analysis, geometry and number
theory, Progress in  Mathematics {\bf 171}, Birkh\"{a}user Boston, Boston, MA, 2000, 29--87.\\
{[D8]} C. Deninger, Number theory and dynamical systems on foliated spaces,
Jahresbericht der  Deutschen Mathematiker-Vereinigung, {\bf 103}, no. 3, 2001, 79--100.\\
{[D9]} C. Deninger, On the nature of the ``explicit formulas" in analytic number theory -- a simple example, 
In: Number Theoretic Methods, Developments in Mathematics, {\bf 8}, Springer, 2002, 97--118.\\
{[D10]} C. Deninger, Arithmetic Geometry and Analysis on Foliated Spaces, arXiv:math/0505354, 2005.\\
{[D11]} C. Deninger, A dynamical systems analogue of Lichtenbaum's conjectures on special values of Hasse-Weil zeta functions, arXiv:math/0605724, 2007.\\
{[D12]} C. Deninger, Analogies between analysis on foliated spaces and
arithmetic geometry, Groups and analysis, London Mathematical Society, Lecture Note
Series,  {\bf 354}, Cambridge Univ. Press, 2008, 174--190.\\
{[D13]} C. Deninger, Dynamical systems for arithmetic schemes, arXiv:1807.06400, v4, 2024.\\
{[E]} A. Erd\'{e}lyi, Higher transcendental functions I, McGraw-Hill, 1953.\\
{[G]} A. Grothendieck, Formulae de Lefschetz et rationalit\'{e} des fonctions $L$, S\'{e}minaire Bourbaki {\bf 279}, 1965.\\
{[KMNT]} J. Kim, M. Morishita, T. Noda, Y. Terashima, On $3$-dimensional foliated dynamical systems and Hilbert type reciprocity law, M\"{u}nster Journal of Mathematics, {\bf 14}, no. 2, 2021, 323--348.\\
{[Ko]} F. Kopei, A foliated analogue of one- and two-dimensional Arakelov
theory, Abhandlungen aus dem Mathematischen Seminar der Universit\"{a}t Hamburg, {\bf 81}, 2011, no. 2, 141--189.\\
{[M]} M. Morishita, Knots and Primes -- An Introduction to Arithmetic Topology, 2nd edition, Universitext, Springer, 2024.\\
{[S]} T. Sunada, Geodesic flows and geodesic random walks, Advanced Studies in Pure Mathematics,  {\bf 3}, 1984, 47--85.\\
\end{flushleft}
{\small
J. A. \'{A}lvarez L\'{o}pez \\
Department of Mathematics and \\
Galician Centre for Math Research and Technology (CITMAga), \\
University of  Santiago de Compostela,\\
15782, Santiago de Compostela, Spain\\
e-mail: jesus.alvarez@usc.es\\
\\
J. Kim:\\
3-18-3, Megurohoncho, Meguro-ku, Tokyo 152-0002, Japan\\
e-mail: res1235@gmail.com\\
\\
M. Morishita:\\
Graduate School of Mathematics, \\
Kyushu University, \\
744, Motooka, Nishi-ku, Fukuoka  819-0395, Japan.\\
e-mail: morishita.masanori.259@m.kyushu-u.ac.jp \\
}
\end{document}